\documentclass[11pt,reqno]{amsart}
\usepackage[T2A]{fontenc}
\usepackage[cp1251]{inputenc}
\usepackage[russian,english]{babel}

\usepackage[pdftex,unicode,colorlinks,linkcolor=blue,
citecolor=red,bookmarksopen,pdfhighlight=/N]{hyperref}

\usepackage{indentfirst}
\usepackage{verbatim}
\usepackage{amsfonts,amssymb,mathrsfs,amscd}
\usepackage{amsmath,latexsym}
\usepackage{amsthm}
\usepackage{verbatim}
\usepackage{eucal}
\usepackage{graphicx}
\usepackage{enumerate}

\theoremstyle{plain} 
\newtheorem{theorem}{Теорема} 
\newtheorem*{theoremA}{Теорема A} 
\newtheorem*{lemmaA}{Лемма A} 
\newtheorem*{thGM}{Теорема Гришина\,--\,Малютиной о малых интервалах}

\newtheorem{propos}{Предложение}[section] 
\newtheorem{corollary}{Следствие}[section] 
\theoremstyle{definition}
 
\newtheorem{remark}{Замечание}[section]

\renewcommand{\leq}{\leqslant} 
\renewcommand{\geq}{\geqslant}
\newcommand{\rad}{\text{\tiny\rm rad}}
\newcommand{\RR}{\mathbb{R}} 
\newcommand{\CC}{\mathbb{C}} 
\newcommand{\NN}{\mathbb{N}}

\DeclareMathOperator{\har}{har} 
\DeclareMathOperator{\mes}{mes}
\DeclareMathOperator{\ess}{ess} 
\DeclareMathOperator{\sbh}{sbh} 
\DeclareMathOperator{\dd}{\,{\mathrm d\!}}

\begin{document}

%        Заголовок статьи.
\title{Одна теорема о малых интервалах\\ для субгармонических функций}
%%\shorttit{\MakeUppercase{Одна теорема  о малых интервалах для субгармонических %%функций}} 
%% верхний колонтитул на нечетных страницах
%% Здесь привести (закомментаренный!) заголовок статьи на английском языке
%% "A Small Intervals Theorem  for Subharmonic Functions''
%%        Автор(ы), инициалы, фамилия
\author{Л.\,А.~Габдрахманова, Б.\,Н.~Хабибуллин}
%%\tauthor{\MakeUppercase{Л.\,А.~Габдрахманова, Б.\,Н.~Хабибуллин}}  %% верхний колонтитул на четных страницах

\date{01/11/2019}

\selectlanguage{russian}
\maketitle

\section{Введение. Основной результат}
\subsection{Предшествующие результаты}\label{Ss1_1} В теории роста целых  и мероморфных  функций на  комплексной плоскости $\CC$ нередко возникает необходимость в оценках интегралов, подынтегральное выражение в которых построено по логарифму модуля таких функций и, возможно,  дополнено некоторой весовой функцией-множителем. Такие оценки относительно множества интегрирования можно  отнести к одному из следующих двух типов:   по различным   интервалам на {\it вещественной оси} $\RR\subset \CC$ и на {\it положительной  
полуоси\/} $\RR^+:=\{x\in \RR\colon x\geq 0\}\subset \RR$
или же    по малым подмножествам на интервалах дуги или луча. Возможна и одновременная  комбинация этих двух типов оценок. 
Среди исходных результатов  первого типа можно рассматривать одну из классических теорем Р.~Неванлинны \cite[стр. 24--27]{RNevanlinna}, изложенную 
в монографии А.\,А.~Гольдберга и И.\,В. Островского  \cite{GO}. 
Истоки второго типа оценок ---   лемма А.~Эдрея и В.~Фукса о малых дугах \cite[{\bf 2}, лемма III, {\bf 9}]{EF}, нашедшая важные применения в теории мероморфных функций, отражённые, в частности, в  \cite[гл. 1, теоремы 7.3, 7.4]{GO}. Вариации на тему леммы о малых дугах ---  лемма о малых интервалах А.\,Ф.~Гришина и М.\,Л.~Содина, доказательство которой, как отмечено авторами в \cite{GrS}, дословно повторяет доказательство леммы Эдрея\,--\,Фукса, и поэтому она приводится в 
\cite{GrS} без доказательства, но с интересным применениями  \cite[лемма 3.2, теорема 3.1]{GrS} к мероморфным функциям конечного порядка, а также к целым функциям  слабо регулярного роста  на луче \cite[\S~3, определение]{GrS}.
Версия этой леммы для субгармонических функций \textit{конечного порядка}  доказана в совместной работе А.\,Ф.~Гришина и Т.\,И.~Малютиной
 \cite{GrM}, посвящённой ряду вопросов теории роста субгармонических функций. Она нашла в  \cite{GrM} неоднократные важные применения в доказательствах ключевых результатов  \cite[теоремы 2, 4]{GrM}. Приведём здесь эту версию в максимально близкой к первоисточнику форме.

\begin{thGM}[{\rm \cite[теорема 8]{GrM}}] Пусть $\boldsymbol \rho \colon (0,+\infty) \to  [0,+\infty)$ --- уточнённый порядок в смысле Валирона  \cite[гл.~III, 6]{Valiron}, \cite[гл.~I, \S~12]{Levin56}, \cite[гл.~II, \S~2]{GO},  т.\,е. функция $\boldsymbol \rho$ дифференцируема  и существуют два  конечных предела 
$\lim_{r\to +\infty}\boldsymbol \rho(r):=\rho\in 0,+\infty)$, 
$\lim_{r\to +\infty}r\boldsymbol \rho'(r)\ln r=0$.  
Пусть $v\not\equiv -\infty$ --- субгармоническая функция формального порядка 
$\boldsymbol \rho$ на $\CC$, т.\,е. существует постоянная $c\in [0,+\infty)$, с которой  
\begin{equation}\label{v}
v(z)\leq c r^{\boldsymbol \rho (r)}\quad \text{при $r=|z|$ для всех  $z\in \CC$},
\end{equation}
$E$ --- измеримое множество на полуоси  $[1,+\infty)$. Тогда существует величина  $ M_{11}$, не зависящая от $R, \theta, E$, такая, что выполняется неравенство
\begin{equation}\label{estGM}
\int_0^R \mathbf{1}_E(t)\bigl|v(te^{i\theta})\bigr| \dd t\leq 
M_{11}\frac{\mes E_R}{R}\ln \frac{4R}{\mes E_R} \,R \cdot R^{\boldsymbol \rho (R)}, \quad E_R:=E\cap [0,R],
\end{equation}
где $\mathbf{1}_E$ --- характеристическая функция подмножества  $E$ в $\RR$, а 
$\mes E_R$ ---  линейная мера Лебега  множества $E_R$. 
 \end{thGM}
В настоящей статье рассматриваются интегралы по подмножествам на интервалах  на луче и  не затрагиваются оценки интегралов по малым дугам на окружностях. В теореме \ref{th1}, сформулированной в подразделе \ref{lo}, мы  распространяем теорему Гришина\,--\,Малютиной о малых интервалах на \textit{субгармонические функции произвольного роста} с полностью контролируемыми постоянными в оценках  с точностью до абсолютных постоянных-множителей. При этом теорема Гришина\,--\,Малютиной усиливается и для субгармонических функций $v$ формального порядка $\boldsymbol \rho $ из \eqref{v}, а именно: аналог постоянной $M_{11}$ из \eqref{estGM} перестает зависеть от субгармонической функции $v$ из \eqref{v} при условии полунормировки  $v(0)\geq 0$, а зависит только от  уточнённого порядка $\boldsymbol \rho $ и постоянной $c$ из правой части \eqref{v}. 
Менее существенное дополнение   в нашем основном результате, --- теореме \ref{th1}, сформулированной в подразделе  \ref{lo}, --- 
допущение произвольной    измеримой функции-множителя $g$ в подынтегральном выражении в интеграле из \eqref{estGM}, хотя иногда это может быть полезным --- см. замечание \ref{remg}.

\subsection{Основные  обозначения и соглашения} 
К этому подразделу можно обращаться по мере необходимости. Наши обозначения могут отличаться от использованных выше. 

 $\RR_{-\infty}:=\{-\infty\}\cup \RR$, $\RR_{+\infty}:=\RR\cup \{+\infty\}$; 
$\RR_{\pm\infty}:=\RR_{-\infty}\cup \RR_{+\infty}$
--- {\it расширенная вещественная ось\/} с естественным отношением порядка, продолженным с $\RR$ до неравенств $-\infty \leq x\leq +\infty$ для любого $x\in \RR_{\pm\infty}$ и порядковой топологией с двумя концами $\pm\infty$;
$S_*:=S\setminus \{0\}$  для $S\subset \RR_{\pm\infty}$, 
в частности,  $\RR_*^+:=(\RR^+)_*=\RR^+\setminus \{0\}$. Определены следующие действия: 
\begin{equation}\label{actR}
\begin{split}
x+(+\infty)=+\infty \text{ для $x\in \RR_{+\infty}$},
\; x+(-&\infty)=-\infty\text{ для }x\in \RR_{-\infty},\; -(\pm\infty)=\mp\infty,
\\
x^+:=\max\{0,x\}=:(-x)^-,\;
x\cdot (\pm\infty):=&\pm\infty=:(-x)\cdot (\mp\infty) \text{ для $x\in \RR_*^+\cup (+\infty)$}, 
\\
\frac{\pm x}{0}:=\pm\infty\text{ для $x\in  \RR_*^+\cup \{+\infty\}$},& \; 
\frac{x}{\pm\infty}:=0\text{ для $x\in  \RR$},  \text{ но $0\cdot (\pm\infty):=0$,}
\end{split}
\end{equation}   
если не оговорено противное. {\it Интервал\/} $I$ на вещественной оси  $\RR$ 
или её расширении $\RR_{\pm\infty}$ ---  {\it связное множество\/}
с {\it концами\/} $\inf I$ и $\sup I$ в $\RR_{\pm\infty}$.
Используем широко распространённые обозначения для {\it открытых интервалов\/}  
$(a,b):=\{x\in \RR_{\pm\infty} \colon a<x<b\}$ 
и  {\it замкнутых интервалов\/} 
$[a,b]:=\{x\in \RR_{\pm\infty} \colon a\leq x\leq b\},$
 а также для  {\it полуоткрытых\/} или {\it полузамкнутых интервалов\/} 
$[a,b):=[a,b]\setminus \{b\}$, $(a,b]:=[a,b]\setminus \{a\}$  
с {\it концами\/} $a,b\in \RR_{\pm\infty}$. 

Для $z\in \CC$  при $r\in \RR^+\cup \{+\infty\}$ через 
$ D(z,r):=\{z' \in \CC \colon |z'-z|<r\}$ обозначаем  {\it открытый круг с центром\/  $z$
радиуса\/ $r$,} $\overline  D(z,r):=\{z' \in \CC \colon |z'-z|\leq r\}$ ---
{\it замкнутый круг с центром\/ $z$ радиуса\/ $r$,} 
$\partial \overline D(z,r):=\{z' \in \CC \colon |z'-z|=r\}$ --- 
{\it окружность  с центром\/ $z$ радиуса\/ $r$.} При $r=0$ из этих определений  $D(z,0)=\varnothing$, $\overline  D(z,0)=\partial \overline D(z,0)=\{z\}$.

 Для $X\subset \RR_{\pm\infty}$ функция $f\colon X\to \RR_{\pm\infty}$ {\it возрастающая\/} (соответственно {\it убывающая\/}) если для любых  $x_1,x_2\in X$ из $x_1\leq x_2$ следует $f(x_1)\leq f(x_2)$ (соответственно $f(x_1)\geq f(x_2)$). 

Для $S\subset \CC$ через $\har(S)$ и  $\sbh(S)$  обозначаем классы всех соответственно {\it гармонических\/} и {\it  субгармонических\/}  функций на каких-либо открытых окрестностях $S$ \cite{Rans}, \cite{HK}. Так, для $z\in \CC$ класс $\har(z):=\har(\{z\})$ состоит из гармонических функций в  окрестностях  точки $z$. 

Через $\lambda$ обозначаем {\it линейную меру Лебега\/} на $\RR$
и её ограничения на {\it $\lambda$-измеримые подмножества\/} $S\subset \RR_{\pm\infty}$.  При этом используем и обозначение $\mes S:=\lambda(S)$, отдавая дань традиции, отражённой выше  в формулировке теоремы Гришина\,--\,Малютиной о малых интервалах.  Понятия \textit{почти всюду,  измеримость\/} и {\it интегрируемость\/} означают  $\lambda$-почти всюду, $\lambda$--измеримость и  $\lambda$--интегрируемость соответственно
Для   \textit{определенной почти всюду\/}   и  \textit{измеримой\/}  функции  $f\colon S\to \RR_{\pm\infty}$ через 
\begin{equation}\label{essg} 
\ess \sup_S f:=\inf \Bigl\{a\in \RR\colon \lambda\bigl(\{x\in S\colon
f(x)>a\}\bigr)=0 \Bigr\}
\end{equation}
обозначаем  {\it существенную  верхнюю грань\/} функции $f$ на $S$. 

Если  $S:=I$ --- интервал с концами $a\leq b$, то для {\it интеграла Лебега\/} от $f$ по интервалу $I$  будем использовать  две равноправные формы обозначений
$$\int_If\dd \lambda:=:\int_a^b f(x) \dd x,$$
если такой интеграл корректно определён. 

Для $r\in \RR^+$ и произвольной функции $v\colon \partial \overline D(0,r)\to \RR_{\pm\infty}$ определим 
\begin{equation}\label{Mvr}
\begin{split}
{\sf M}_v(r)&:=\sup_{|z|=r}v(z);\\
{\sf C}_v(r)&:=\frac{1}{2\pi}\int_0^{2\pi} v(re^{is})\dd s 
\end{split}
\end{equation}
--- {\it среднее по окружности\/ $\partial \overline{D}(0,r)$ функции\/} $v$ при условии, что интеграл Лебега в правой части  корректно определён. В частности, 
${\sf M}_v(0)={\sf C}_v(0)=v(0)$. Свойства характеристик 
${\sf M}_v$ и  ${\sf C}_v$ в случае субгармонической функции $v$ см. в 
\cite[2.6]{Rans}, \cite[2.7]{HK}

\subsection{Формулировка основного результата}\label{lo}

\begin{theorem}[{о малых интервалах}]\label{th1} 
Существует абсолютная постоянная $a\geq 1$, с которой  
\begin{enumerate}[{\rm (i)}]
\item[{\rm [{\bf u}]}]\label{li} для любой  субгармонической на $\CC$ функции $u\not\equiv -\infty$, 
\item[{\rm [{\bf r}]}]\label{lii} для любых чисел $0\leq r_0\leq  r<R<+\infty$,
\item[{\rm [{\bf E}]}]\label{liii} для любого измеримого подмножества $E\subset [r,R]$,
\item[{\rm [{\bf g}]}]\label{liv} для любой измеримой функции $g$
на $E$, 
\item[{\rm [{\bf b}]}]\label{lb} для любого числа  $b\in (0,1]$
\end{enumerate} 
имеет место неравенство 
\begin{equation}\label{iend}
\begin{split}
\int_{E}{\sf M}_{|u|} g\dd \lambda&\leq \left(\frac{a}{b}\ln \frac{a}{b}\right)
 \Bigl( {\sf M}_u\bigl((1+b)R\bigr)-2\min \bigl\{0, {\sf C}_u(r_0)\bigr\}\Bigr)
 \Bigl(\ess \sup_{E} |g|\Bigr)\\
&\times \left(\mes E+\min\{\mes E, 3bR\} \ln \frac{3ebR}{\min\{\mes E, 3bR\}}\right)
\end{split}
\end{equation}
\end{theorem}
\begin{remark}\label{remes}
Последнюю скобку-множитель 
\begin{equation}\label{M}
M\overset{\eqref{iend}}{:=}\mes E+\min\{\mes E, 3bR\} \ln \frac{3ebR}{\min\{\mes E, 3bR\}}
\end{equation}
часто можно <<сжать>>. Так, 
 \begin{equation}\label{E<R}
M\overset{\eqref{M}}{\leq} 2 \mes E \ln \frac{3ebR}{\mes E},
\quad\text{eсли $\mes E\leq 3bR$,}
\end{equation}
а в противном случае 
\begin{equation}\label{E>R}
M\overset{\eqref{M}}{=}\mes E+ 3bR\leq 2\mes E
\quad\text{при $\mes E> 3bR$,}
\end{equation}
где как в случае  \eqref{E<R}, так и в случае \eqref{E>R} число $2$ можно перекинуть как множитель в абсолютную постоянную $a$. 
Кроме того,  всегда $\mes E\leq R$, поскольку $E\subset [r,R]$. При этом 
\begin{equation}\label{x}
(x,t)\longmapsto x\ln \frac{t}{x} 
\text{\it  --- возрастающая функция по\/ $t\in \RR_*^+$ и по\/  $x\in [0,t/e]$.}  
\end{equation}
 Исходя из этого, \textit{при любом\/} $b\in (0,1]$ имеем
\begin{multline}\label{Ml}
M \overset{\eqref{M}}{\leq} 
\mes E+\min\{\mes E, 3bR\} \ln \frac{\max\{3ebR,eR\}}{\min\{\mes E, 3bR\}}
\overset{\eqref{x}}{\leq} \mes E+\mes E \ln \frac{3eR}{\mes E}\\
=(1+\ln 3) \mes E+ \mes E\ln \frac{eR}{\mes E}\leq (2+\ln 3) \mes E\ln \frac{eR}{\mes E},
\end{multline}
где множитель  $2+\ln 3$ можно перебросить в абсолютную постоянную $a$. Конечно же, замена значения $M$ из правой части \eqref{M}  на правую часть из \eqref{Ml}
в оценке сверху \eqref{iend} может вести к значительному огрублению оценки  \eqref{iend} при малых значениях $b>0$.
\end{remark}

\begin{remark}\label{remg}
Если функция $g\geq 0$ в {[{\bf g}]} \textit{возрастающая,\/} то в неравенстве \eqref{iend} 
величину $\ess \sup_{E} |g|$ можно заменить на $g(R)$, а если \textit{убывающая,\/} то на $g(r)$.

\end{remark}
Покажем, что теорема Гришина\,--\,Малютиной 
содержится в теореме \ref{th1}. 

\begin{proof}[Вывод теоремы Гришина\,--\,Малютиной]  По условию  $E\subset [1,+\infty)$. Положим $g\equiv 1$ на $\RR^+$, а $b:=1$. По теореме \ref{th1} о малых интервалах с учётом замечания \ref{remes} с оценкой \eqref{Ml} найдётся абсолютная постоянная $A\geq 1$, с которой оценку \eqref{iend}  при $r_0=r=1<R$  для субгармонической функции $v\not\equiv -\infty$, удовлетворяющей  \eqref{v}, можно записать в виде 
\begin{multline}\label{GKh}
\int_0^R \mathbf{1}_E(t)\bigl|v(te^{i\theta})\bigr| \dd t \leq 
\int_{E_R}{\sf M}_{|v|} \dd \lambda \\
\leq A\Bigl( c(2R)^{\boldsymbol \rho (2R)}-2\min\bigl\{ 0,{\sf C}_v(1)\bigr\}\Bigr)
\mes E_R\ln \frac{eR}{\mes E_R}.
\end{multline}
Из свойств уточнённого порядка $\mathbb \rho$ следует существование 
постоянной $C_{\boldsymbol \rho}$, с которой 
$(2R)^{\boldsymbol \rho (2R)}\leq C_{\boldsymbol \rho}R^{\boldsymbol \rho (R)}$ для всех $R\geq 1$ \cite[гл.~II, теорема 2.2, следствие]{GO}, \cite[гл.~1, \S~12, 2)]{Levin56}. Таким образом, из \eqref{GKh}
получаем
\begin{equation*}
\int_0^R \mathbf{1}_E(t)\bigl|v(te^{i\theta})\bigr| \dd t\leq 
A\left(cC_{\boldsymbol \rho}-2\frac{\min\bigl\{ 0,{\sf C}_v(1)\}}{\inf_{r\geq 1}
r^{\boldsymbol \rho (r)}}\right)\mes E_R\ln \frac{eR}{\mes E_R}  \cdot R^{\boldsymbol \rho (R)}, 
\end{equation*}
откуда следует \eqref{estGM}.
\end{proof}

\section{Доказательство теоремы \ref{th1}} 
\setcounter{equation}{0} 

\subsection{Оценка модуля субгармонической функции}\label{Intr}
Нам потребуется одна из оценок снизу для субгармонических функций, установленная  в \cite{Kh84d}, \cite{Kh84} в  гораздо более общей форме, чем необходимо 
для наших целей.
\begin{theoremA}[{\rm \cite[теорема 1.2]{Kh84d}, \cite[теорема 2]{Kh84}}]\label{thA}
Существует такая абсолютная постоянная $a_0\geq 1$, что для любой убывающей непрерывной функции  $a\colon \RR^+\to \RR_*^+$ и для любой возрастающей непрерывной функции   $A\colon \RR^+\to \RR_*^+$, связанных соотношениями 
$1\leq aA\leq A$ на $\RR^+$,  для каждой  функции 
\begin{equation}\label{uh0}
u\in \sbh(\CC)\cap \har (0) \quad\text{со значением } u(0)=0 
\end{equation}
 найдётся не более чем счётное множество  кругов 
\begin{equation}\label{Dj}
\{D(z_j,t_j)\}_{j\in N}, \quad N\subset \NN:=\{1, 2, \dots \},
\end{equation} 
для которого  неравенства 
$$
u(z)\geq -\frac{a_0}{a(2|z|)} {\sf M}_u\bigl((1+a(|z|))|z|\bigr)\ln \bigl(a_0A(|z|)\bigr)
$$
выполнены при всех $z\in \CC\setminus \bigcup_{j\in N}D(z_j,t_j)$ и при этом
$$
\sum_{ D(0,r)\cap D(z_j,t_j)\neq \varnothing} t_j \leq \frac{1}{a(2r)} \int_0^r\frac{1}{A(t)} \quad\text{для любых $r\in \RR_*^+$}.
$$
\end{theoremA}

\begin{corollary}\label{prsn}
Существует такая абсолютная постоянная $a_0\geq 1$, что 
для любой субгармонической функции  $u$, удовлетворяющей  ограничениям из \eqref{uh0},  и для любого числа $B>1$ найдётся не более чем счётное множество кругов \eqref{Dj}, для которого
\begin{equation}\label{esl}
u(z)\overset{\eqref{Mvr}}{\geq} -\bigl(a_0B\ln (a_0B)\bigr) {\sf M}_u\bigl((1+1/B)|z|\bigr), 
\end{equation}
при всех $z\in \CC\setminus \bigcup_{j\in N}D(z_j,t_j)$ и при этом
\begin{equation}\label{eslb}
\sum_{D(r)\cap D(z_j,t_j)\neq \varnothing}t_j\leq \frac{1}{B}\, r\quad
 \text{для любых  $r\in \RR_*^+$.}
\end{equation}
\end{corollary}
\begin{proof}
Предложение \ref{prsn} сразу следует из   теоремы A  
%\cite[теорема 1.2]{Kh84d}, \cite[теорема 2]{Kh84} 
при  постоянных функциях $A(t)\equiv B^2$ и $a(t)\equiv 1/B$ при всех  $t\in \RR^+$. 
\end{proof} 
\begin{propos}\label{disk}
Для любого числа  $b\in (0,1]$ при 
\begin{equation}\label{B}
B>2+\frac{2}{b}
\end{equation}
для каждого множества кругов 
\eqref{Dj} с ограничением \eqref{eslb}, для всякой  точки $z\in \CC\setminus\{0\}$ найдётся число $r_z>0$, удовлетворяющее условиям 
\begin{equation}\label{rqabD}
0< r_z < b|z|,
\quad  \biggl(\; \bigcup_{j\in N}D(z_j,r_j)\biggr)
\bigcap \partial D(z,r_z)
=\varnothing,
\end{equation}
т.\,е. окружность $\partial \overline D(z,r_z)$ при некотором $r_z>0$ содержится в  круге $D\bigl(z, b|z|\bigr)$ и не пересекается ни с одним из кругов $D(z_j,r_j)$ с $j\in N$.
\end{propos} 
Предложение \ref{disk} элементарно, но в целях контроля за соотношениями между постоянными $b$ и $B$ всё же приведём его  стандартное 
\begin{proof}
Достаточно рассмотреть положительную точку $z\in \RR_*^+$. 
При   ограничении  \eqref{eslb} и значении $B$ из \eqref{B} имеем
\begin{equation*}
 \sum_{D(z,bz)\cap D(z_j,t_j)\neq \varnothing}2t_j
\leq \sum_{D((1+b)z)\cap D(z_j,t_j)\neq \varnothing}2t_j
\overset{ \eqref{eslb}}{\leq} \frac{2}{B}\, (1+b)z
\overset{\eqref{B}}{<} b z.
\end{equation*}
Отсюда для всех кругов $D(z_j,t_j)$ с $D(z,bz)\cap D(z_j,t_j)\neq \varnothing$,
сумма длин интервалов 
\begin{equation}\label{prdb}
\Bigl(z+\bigl(D(z_j,t_j)-z\bigr)e^{-i\arg (z_j-z)}\Bigr)\bigcap \bigl[0,bz\bigr), 
\end{equation}
являющихся  круговыми проекциями с центром $z$ в $D(z,bz)$ диаметров кругов $D(z_j,t_j)$ на радиус $[z,bz)\subset \RR_*^+$ круга $D(z,bz)$, меньше, чем длина $bz$ радиуса круга $D(z,bz)$. Следовательно, на радиусе $[z,bz)\subset \RR_*^+$ есть точка $z+r_z$, не принадлежащая ни одному из интервалов  \eqref{prdb}. Это сразу влечёт за собой  \eqref{rqabD}.
\end{proof}

\begin{propos}\label{lemrep}
Существует такая  абсолютной постоянной $A_0\geq 1$, что 
\begin{enumerate}
\item[{\rm [u]}] для любой субгармонической функции $u$ из \eqref{uh0} 
c мерой Рисса $\mu_u:=\frac{1}{2\pi}\Delta u$, где $\Delta$ --- оператор Лапласа, действующий в смысле теории обобщённых функций,
\item[{\rm [{\bf b}]}]  для любого числа  $b\in (0,1]\subset \RR^+$
\end{enumerate} 
имеет место неравенство  
\begin{equation}\label{uCMl}
{\sf M}_{|u|}(r)\leq \int_{(1-b)r}^{(1+b)r}
\ln\frac{b r}{|t-r|} \dd \mu_u^{\rad}(t)
+\Bigl(\frac{A_0}{b}\ln \frac{A_0}{b}\Bigr){\sf M}_u\bigl((1+2b)r\bigr)\quad \text{при всех $r\in \RR^+$,} 
\end{equation}
где $\mu^{\rad}(t):=\mu\bigl(\overline D(0,t)\bigr)$, $t\in \RR^+$,
--- считающая функция положительной меры $\mu$ на $\CC$.
\end{propos}
\begin{proof} По определению \eqref{Mvr} из принципа максимума 
\begin{equation}\label{est_ab}
u(z)\leq {\sf M}_u\bigl(|z|\bigr)\leq {\sf M}_u\bigl((1+2b)|z|\bigr) \quad\text{при всех } z\in \CC.
\end{equation} 
Для абсолютной постоянной $a_0\geq 1$ из следствия  \ref{prsn} положим 
\begin{equation}\label{BC}
B\overset{\eqref{B}}{:=}3+\frac{3}{b}, \quad 
1\leq c_b:=a_0B\ln (a_0B)=\frac{3a_0(b+1)}{b} \ln \frac{3a_0(b+1)}{b}\leq\frac{6a_0}{b} \ln \frac{6a_0}{b} . 
\end{equation}
Пусть $z\in \CC\setminus \{0\}$. Исходя из оценки снизу \eqref{esl} следствия  \ref{prsn}
для функции $u$ и предложения  \ref{disk},  найдется такое число $r_z\in \bigl(0, b|z|\bigr)$, удовлетворяющее  \eqref{rqabD}, что имеет место  оценка снизу 
\begin{multline}\label{es_rzM}
u(z')\geq - {\sf M}_u\bigl((1+1/B)|z'|\bigr)a_0B \ln (a_0B)
\overset{\eqref{BC}}{\geq} -c_b {\sf M}_u\bigl((1+1/B)(1+b)|z|\bigr)\\
\overset{\eqref{BC}}{\geq} -c_b {\sf M}_u\bigl((1+2b)|z|\bigr)
\quad \text{для всех $z'\in \partial\overline D(z,r_z)$}.
\end{multline}
По формуле Пуассона\,--\,Йенсена \cite[3.7]{HK}, \cite[4.5]{Rans} для круга $D(z,r_z)$ имеем
\begin{equation}\label{fPJv}
u(z)=-\int_{D(z,r_z)} \ln\Bigl|\frac{r_z^2}{r_z(z-z')}\Bigr|\dd \mu_u(z')+h(z)=:G(z)+h(z),
\end{equation}
где наименьшая гармоническая мажоранта $h$ функции $v$ в круге $D(z,r_z)$ ввиду \eqref{est_ab} и \eqref{es_rzM} удовлетворяет оценке 
\begin{equation}\label{esh}
\bigl|h(z)\bigr|\leq c_b {\sf M}_u\bigl((1+2b)|z|\bigr).
\end{equation}
Поскольку функция Грина для круга под знаком интеграла c полюсом $z$ в центре круга в промежуточном равенстве  из \eqref{fPJv} положительна, имеем 
\begin{multline*}
\bigl|G(z)\bigr|= \int_{D(z,r_z)} \ln\frac{r_z}{|z-z'|}\dd \mu_u(z')
\leq \int_{D(z,r_z)}\ln\frac{b|z|}{|z-z'|}\dd \mu_u(z') 
\leq \int_{D(z,b|z|)}\ln\frac{b|z|}{\bigl||z'|-|z|\bigr|}
\dd \mu_u(z')
\\
\leq \int\limits_{(1-b)|z|\leq |z'|<(1+b)|z|}\ln\frac{b|z|}{\bigl||z'|-|z|\bigr|}
\dd \mu_u(z')
=\int\limits_{(1-b)|z|}^{(1+b)|z|}
\ln\frac{b|z|}{\bigl|t-|z|\bigr|} \dd \mu_u^{\rad}(t).
\end{multline*} 
Из последней оценки и из \eqref{esh} ввиду \eqref{fPJv} получаем 
\begin{equation*}
\bigl|u(z)\bigr|\leq \bigl|G(z)\bigr|+\bigl|h(z)\bigr|\leq 
\int\limits_{(1-b)|z|}^{(1+b)|z|}
\ln\frac{b|z|}{\bigl|t-|z|\bigr|} \dd \mu_u^{\rad}(t)+
c_b {\sf M}_u\bigl((1+2b)|z|\bigr)
\end{equation*}
Здесь правая часть не зависит от аргумента $\arg z$,
что при $r:=|z|>0$ по \eqref{Mvr} даёт  \eqref{uCMl} с $A_0\overset{\eqref{BC}}{:=}6a_0$ . При $r=0$ неравенство 
\eqref{uCMl} очевидно, поскольку $u(0)=0$.
\end{proof}

\subsection{Интегральные оценки}\label{1stl}
Если $\mes E=0$, то в рамках соглашений \eqref{actR} о действиях в $\RR_{\pm\infty}$
обе части неравенства \eqref{iend} равны нулю и оно верно. Поэтому далее $\mes E>0$.
При этом если $\ess \sup_{E} |g|=+\infty$, то согласно  \eqref{actR}
правая часть  неравенства \eqref{iend} равна $+\infty$ и неравенство \eqref{iend} верно.
В противном  случае сразу можем <<убрать>> функцию $g$ из подынтегрального выражения в интеграле из  левой части \eqref{iend}:
\begin{equation}\label{Mg}
\int_E {\sf M}_{|u|}g \dd \lambda
\leq  
\int_r^R {\sf M}_{|u|}(x) {\boldsymbol 1_{E}} \dd x
 \cdot\ess \sup_{E} |g|.
\end{equation}
Основные этапы доказательства содержит случай 
\subsubsection{\bf Функция $u$ из \eqref{uh0}}\label{sssu} 
Пока считаем, что $b\in (0,1/2]$. 
По предложению \ref{lemrep} для любой функции $u$ из \eqref{uh0} для интеграла из правой части \eqref{Mg} имеет место неравенство 
\begin{multline}\label{intrR}
\int_r^R {\sf M}_{|u|}(x) {\boldsymbol 1_{E}}(x) \dd x
\leq 
\int_r^R\left(\int_{(1-b)x}^{(1+b)x}
\ln\frac{b x}{|t-x|} \dd \mu_u^{\rad}(t)
+\Bigl(\frac{A_0}{b}\ln \frac{A_0}{b}\Bigr){\sf M}_u\bigl((1+2b)x\bigr)\right) {\boldsymbol 1_{E}}(x) \dd x\\
\leq \int_r^R\int_{(1-b)x}^{(1+b)x}
\ln\frac{b x}{|t-x|} \dd \mu_u^{\rad}(t){\boldsymbol 1_{E}}(x) \dd x
+\Bigl(\frac{A_0}{b}\ln \frac{A_0}{b}\Bigr){\sf M}_u\bigl((1+2b)R\bigr) \mes E.
 \end{multline}
 Обозначим повторный интеграл в правой части \eqref{intrR} 
как 
\begin{equation}\label{I}
I:= \int_r^R\int_{(1-b)x}^{(1+b)x}
\ln\frac{b x}{|t-x|} \dd \mu_u^{\rad}(t){\boldsymbol 1_{E}}(x) \dd x
=\iint_{S}  \boldsymbol{1}_E(x) \ln^+\frac{bx}{|t-x|}  \dd \mu_u^{\rad}(t)\dd x 
\end{equation}
где $\ln^+:=\max\{0,\ln\}$ и  двойной интеграл в правой части берётся по  множеству 
\begin{equation*}
S:=\Bigl\{(x,t)\in \RR^2\colon x\in [r,R], (1-b)x\leq t\leq(1+b)x  \Bigr\}
\subset \RR^2, 
\end{equation*}
содержащемуся  при $b\in (0,1/2]$ в множестве
\begin{equation*}
T:=\Bigl\{(x,t)\in \RR^2\colon t\in \bigl[(1-b)r,(1+b)R\bigr], (1-2b)t\leq x\leq (1+2b)t  \Bigr\}.
\end{equation*} 
Отсюда  согласно \eqref{I} при $b\in (0,1/2]$ получаем 
\begin{multline}\label{Iest}
I \leq \iint_{T} \boldsymbol{1}_E(x) \ln^+\frac{bx}{|t-x|} 
 \dd \mu_u^{\rad}(t)\dd x
=\int_{(1-b)r}^{(1+b)R}\int_{(1-2b)t}^{(1+2b)t}
 \boldsymbol{1}_E(x) \ln^+\frac{bx}{|t-x|}  
\dd x \dd \mu_u^{\rad}(t)
\\
\leq \Bigl(\mu_u^{\rad}\bigl((1+b)R\bigr)-\mu_u^{\rad}\bigl((1-b)r\bigr)
\Bigr)
\sup_{(1-b)r\leq t\leq (1+b)R}\int_{(1-2b)t}^{(1+2b)t}
 \boldsymbol{1}_E(x) \ln^+\frac{bx}{|t-x|}  \dd x\\
\leq  \mu_u^{\rad}\bigl((1+b)R\bigr)
\sup_{(1-b)r\leq t\leq (1+b)R} \int_0^{(1+b)(1+2b)R}\boldsymbol{1}_E(x)
 \ln^+\frac{b(1+b)(1+2b)R}{|t-x|}  \dd x\\
\leq  \mu_u^{\rad}\bigl((1+b)R\bigr)
\sup_{ t\in \RR}\int_{-\infty}^{+\infty}\boldsymbol{1}_E(x)
 \ln^+\frac{3bR}{|x-t|}  \dd x\\=\mu_u^{\rad}\bigl((1+b)R\bigr)
\sup_{ t\in \RR}\int_{-\infty}^{+\infty}\boldsymbol{1}_{(E-t)}(s)
 \ln^+\frac{3bR}{|s|}  \dd s,
\end{multline}
где множитель перед $\sup$
стандартно оценивается  сверху при $b\in (0,1/2]$ как 
\begin{multline}\label{murR}
\mu_u^{\rad}\bigl((1+b)R\bigr)\leq\frac{1}{\ln\frac{1+2b}{1+b}}
\int_{(1+b)R}^{(1+2b)R}\frac{\mu_u^{\rad}(t)}{t}\dd t
\leq \frac{4(1+b)}{b} \int_0^{(1+2b)R}\frac{\mu_u^{\rad}(t)}{t}\dd t\\
= \frac{6}{b} {\sf C}_u\bigl((1+2b)R\bigr)
\leq \frac{6}{b}{\sf M}_u\bigl(1+2b)R\bigr),
\end{multline}
поскольку в этом п.~\ref{sssu} функция $u$ удовлетворяет требованиям из  \eqref{uh0}. 
Для оценки последнего интеграла в \eqref{Iest} будет использована 
\begin{lemmaA}[{\rm \cite[лемма 7.2]{GO}}] Пусть 
$f\colon I\to \RR_{\pm\infty}$ --- чётная интегрируемая функция на  симметричном относительно нуля открытом  интервале $I\subset \RR$,
убывающая на интервале $I\cap \RR_*^+$, $E\subset I$ --- измеримое подмножество. Тогда 
\begin{equation}\label{intI}
\int_E f\dd \lambda\leq 
2\int_0^{(\mes E)/2} f(x)\dd x.
\end{equation}
\end{lemmaA}
Применяя лемму  A  с  
$f\colon x\mapsto \ln^+ \frac{3bR}{|x|}$,  $x\in \RR$, и множеством 
$E-t$ в роли $E$
к оценке сверху интегралов из правой части \eqref{Iest}, получим  
\begin{multline}\label{estit}
\int_{-\infty}^{+\infty}\boldsymbol{1}_{(E-t)}(s)
 \ln^+\frac{3bR}{|s|}  \dd s 
=\int_{(E-t)} \ln^+\frac{3bR}{|s|}  \dd \lambda(s)
\leq 2\int_0^{\mes (E-t)/2} \ln^+ \frac{3bR}{x} \dd x
\\=2\int_0^{\mes E/2} \ln^+ \frac{3bR}{x} \dd x
=2\int_0^{\min\{\mes E/2, 3bR\}} \ln^+ \frac{3bR}{x} \dd x
\\
=\min\{\mes E, 6bR\} \ln \frac{6ebR}{\min\{\mes E, 6bR\}},
\end{multline}
где последнее выражение не зависит от переменной $t$. Таким образом, последовательно применяя оценку \eqref{intrR}, обозначение \eqref{I} для интеграла $I$ с оценками сверху\eqref{Iest} и \eqref{murR},  а также  \eqref{estit}, в случае функции $u$ из \eqref{uh0} при $b\in (0,1/2]$ получаем   
\begin{multline*}
\int_{E}{\sf M}_{|u|} g\dd \lambda\leq 
\max\left\{\Bigl(\frac{A_0}{b}\ln \frac{A_0}{b}\Bigr), \frac{6}{b}\right\}{\sf M}_u\bigl((1+2b)R\bigr) \ess \sup_{E} |g| \\
\times \left(\mes E+\min\{\mes E, 6bR\} \ln \frac{6ebR}{\min\{\mes E, 6bR\}}\right)
\end{multline*}
Если здесь вместо числа $b\in (0,1/2]$ рассмотреть число $2b\in (0,1]$, сохранив за ним прежнее обозначение $b$, а также положить 
$a:=2(A_0+6)$, то получим в точности неравенство  \eqref{iend}, поскольку для субгармонической функции $u$ из  \eqref{uh0} в силу $u(0)=0$ имеем ${\sf C}_u(r_0)\geq 0$ и $\min \bigl\{0, {\sf C}_u(r_0)\bigr\}=0$.  Осталось рассмотреть 
\subsubsection{\bf Случай произвольной субгармонической функции $u\not\equiv -\infty$}\label{uinfty} При $u(0)=-\infty$ и $r_0=0$ имеем 
$\min \bigl\{0, {\sf C}_u(r_0)\bigr\}=-\infty$. Следовательно, правая часть неравенства  
\eqref{iend} в рамках правил действий \eqref{actR} равна $+\infty$ и неравенство \eqref{iend} верно. 

Пусть $r_0>0$. 
Тогда от произвольной  субгармонической функции $u\not\equiv -\infty$
 можно с помощью её гармонического продолжения внутрь круга $D(0,r_0)$   и вычитания из этого продолжения среднего по окружности 
$c_0\overset{\eqref{Mvr}}{:=}{\sf C}_u(r_0)$ перейти к субгармонической функции 
$u_0$, равной $u-c_0$ на $\CC\setminus D(0,r_0)$, гармонической в круге $D(0,r_0)$ и равной нулю в нуле. Таким образом, эта преобразованная функция $u_0$ удовлетворяет требованиям из \eqref{uh0} и для нее справедливо неравенство  \eqref{iend} с заменой $u$ на $u_0$. Таким образом, применяя \eqref{iend} к $u_0$, с теми же скобками, за исключением второй, в правой  части имеем
\begin{multline*}
\int_{E}{\sf M}_{|u|} |g|\dd \lambda\leq 
\int_{E}{\sf M}_{|u_0|} |g|\dd \lambda
+\int_{E}{\sf M}_{|c_0|} |g|\dd \lambda\\
\leq \Bigl(\dots\Bigr)\Bigl( {\sf M}_u\bigl((1+b)R\bigr)-c_0\Bigr)
  \Bigl(\dots\Bigr) \times \Bigl(\dots\Bigr)+
\int_{E}|c_0| |g|\dd \lambda\\
\leq \Bigl(\dots\Bigr)\Bigl( {\sf M}_u\bigl((1+b)R\bigr)-\bigl(c_0-|c_0|\bigr)\Bigr)
  \Bigl(\dots\Bigr) \times \Bigl(\dots\Bigr)
\end{multline*}
где $c_0-|c_0|=-2c_0^-=2\min \bigl\{0,{\sf C}_u(r_0)\bigr\}$. Это даёт неравенство 
\eqref{iend} для $r_0>0$. Наконец, к случаю $u(0)>-\infty$ и $r_0=0$ можно перейти от последнего разобранного случая, устремляя $r_0>0$ к нулю, с  $\min \bigl\{0,{\sf C}_u(r_0)\bigr\} \to \min  \bigl\{0,u(0)\bigr\}>-\infty$. 

Теорема \ref{th1} доказана.

\begin{remark}
В случае, когда $u=\ln |f|$, где $f$ --- целая ненулевая функция с $|f(0)|=0$, в п.~\ref{uinfty} можно,  конечно, поступить и иначе, --- более традиционно, --- а именно: вместо функции $f$ рассмотреть результат деления функции $f$ на одночлен $z\mapsto \frac{f^{(k)}(0)}{k!}z^k$, где $k\geq 1$ --- кратность корня функции $f$ в нуле.  Но тогда вместо дополнительного слагаемого $-2\min \bigl\{0,{\sf C}_u(r_0)\bigr\}$ в правую часть придётся приписать слагаемое  порядка роста $k\ln R+O(1)$
что при $k\geq 1$ и росте $R\to +\infty$, очевидно, строго грубее, чем постоянная 
$-2\min \bigl\{0,{\sf C}_u(r_0)\bigr\}$.
\end{remark}

\renewcommand{\refname}{Литература}

\end{document}